\documentclass[12pt]{article}
\usepackage{amssymb}
\usepackage{amsmath}

\begin{document}

\begin{center}
\textbf{\large The Bohr--P\'al Theorem and the Sobolev Space
$W_2^{1/2}$}
\end{center}

\begin{center}
Vladimir Lebedev
\end{center}

\begin{quotation}
{\small \textbf{Abstract.} The well-known Bohr--P\'al theorem
asserts that for every continuous real-valued function $f$ on
the circle $\mathbb T$ there exists a change of variable, i.e.,
a homeomorphism $h$ of $\mathbb T$ onto itself, such that the
Fourier series of the superposition $f\circ h$ converges
uniformly. Subsequent improvements of this result imply that
actually there exists a homeomorphism that brings $f$ into the
Sobolev space $W_2^{1/2}(\mathbb T)$. This refined version of
the Bohr--P\'al theorem does not extend to complex-valued
functions. We show that if $\alpha<1/2$, then there exists a
complex-valued $f$ that satisfies the Lipschitz condition of
order $\alpha$ and at the same time has the property that
$f\circ h\notin W_2^{1/2}(\mathbb T)$ for every homeomorphism
$h$ of $\mathbb T$.

2010 \emph{Mathematics Subject Classification}: 42A16.

\emph{Key words and phrases}: harmonic analysis, homeomorphisms
of the circle, superposition operators, Sobolev spaces.}
\end{quotation}

\quad

\begin{center}
\textbf{1. Introduction}
\end{center}

For an arbitrary integrable function $f$ on the circle $\mathbb
T=\mathbb R/2\pi\mathbb Z$ (where $\mathbb R$ is the real line
and $\mathbb Z$ is the group of integers) consider its Fourier
series:
$$
f(t)\sim\sum_{k\in\mathbb Z}\widehat{f}(k)e^{ikt}, \qquad t\in\mathbb T.
$$
Recall that the Sobolev space $W_2^{1/2}(\mathbb T)$ is the
space of all (integrable) functions $f$ with
$$
\sum_{k\in\mathbb Z}|\widehat{f}(k)|^2
|k|<\infty.
$$
Let $C(\mathbb T)$ be the space of all continuous functions on
$\mathbb T$.

The well-known Bohr--P\'al theorem states that for every
real-valued function $f\in C(\mathbb T)$ there exists a
homeomorphism $h$ of the circle $\mathbb T$ onto itself, such
that the superposition $f\circ h$ belongs to the space
$U(\mathbb T)$ of uniformly convergent Fourier series. (The
theorem was obtained in a somewhat weaker form by J. P\'al in
[11], and in the final form by H. Bohr in [2].) The original
method of proof of this result uses conformal mappings and in
fact allows (see [9, Sec. 3]) to obtain the following
representation:
$$
f\circ h= g+\psi, \qquad g\in W_2^{1/2}\cap C(\mathbb T), \quad \psi\in V\cap C(\mathbb T),
\eqno(1)
$$
where $V(\mathbb T)$ is the space of functions of bounded
variation on $\mathbb T$. It is well-known that both
$W_2^{1/2}\cap C(\mathbb T)$ and $V\cap C(\mathbb T)$ are
subsets of $U(\mathbb T)$, thus (1) implies $f\circ h\in
U(\mathbb T)$.

A substantial improvement of the Bohr--P\'al theorem was
obtained by A. A. Sahakian [12, Corollary 1], who showed that
if $a(n), n=0, 1, 2, \ldots,$ is a given positive sequence
satisfying the condition $\sum_n a(n)=\infty$ and a certain
condition of regularity, then for every real-valued $f\in
C(\mathbb T)$ there is a homeomorphism $h$ such that
$\widehat{f\circ h}(k)=O(a(|k|))$. An immediate consequence of
Sahakian's result is that the term $\psi$ in (1) can be
omitted, i.e., the following refined version of the Bohr--P\'al
theorem holds: for every real-valued $f\in C(\mathbb T)$ there
exists a homeomorphism $h$ of $\mathbb T$ onto itself, such
that $f\circ h\in W_2^{1/2}(\mathbb T)$. This refined version
also follows from a result on conjugate functions, obtained by
W. Jurkat and D. Waterman in [4] (see also [3, Theorem 9.5]).
We note that Sahakian's result is obtained by purely real
analysis technique whereas Jurkat and Waterman use an approach
similar to the one used by Bohr and P\'al. A very short proof
of the refined version of the Bohr--P\'al theorem was
communicated to the author by A. Olevski\v{\i}, see [7, Sec.
3].

Another improvement of the Bohr--P\'al theorem was obtained by
J.-P. Kahane and  Y. Katznelson [6] (see also [9], [5]). These
authors showed that if $K$ is a compact family of functions in
$C(\mathbb T)$, then there exists a homeomorphism $h$ of
$\mathbb T$ such that $f\circ h\in U(\mathbb T)$ for all $f\in
K$. This result naturally leads to a question if it is possible
to attain the condition $f\circ h\in W_2^{1/2}(\mathbb T)$ for
all $f\in K$. This question was posed by A. Olevski\v{\i} in
[10]. A negative answer was obtained by the author of this work
in [7, Theorem 4], it turns out that, given a real-valued $u\in
C(\mathbb T)$, the property that for every real-valued $v\in
C(\mathbb T)$ there is a homeomorphism $h$ such that both
$u\circ h$ and $v\circ h$ are in $W_2^{1/2}(\mathbb T)$ is
equivalent to the boundness of variation of $u$. Thus, in
general, there is no single change of variable which will bring
two real-valued functions in $C(\mathbb T)$ into
$W_2^{1/2}(\mathbb T)$. Certainly this amounts to the existence
of a complex-valued $f\in C(\mathbb T)$ such that $f\circ
h\notin W_2^{1/2}(\mathbb T)$ for every homeomorphism $h$ of
$\mathbb T$.

The purpose of this work is to show that there exists a
complex-valued function $f$ that is \emph{very smooth} but at
the same time has the property that $f\circ h\notin
W_2^{1/2}(\mathbb T)$ for every homeomorphism $h$ of $\mathbb
T$.

Note that, as one can easily verify (see, e.g., [7], Sec. 3),
the following two semi-norms
$$
\|f\|_{W_2^{1/2}(\mathbb T)}=\bigg(\sum_{k\in\mathbb Z}|\widehat{f}(k)|^2 |k|\bigg)^{1/2},
$$
$$
|\|f|\|_{W_2^{1/2}(\mathbb T)}=\bigg(\int_0^{2\pi}\frac{1}{\theta^2}
\int_0^{2\pi} |f(t+\theta)-f(t)|^2 dt d\theta\bigg)^{1/2}
\eqno(2)
$$
are equivalent semi-norms on $W_2^{1/2}(\mathbb T)$, i.e., $f$
is in $W_2^{1/2}(\mathbb T)$ if and only if
$|\|f|\|_{W_2^{1/2}(\mathbb T)}<\infty$, and
$c_1\|f\|_{W_2^{1/2}(\mathbb T)}\leq|\|f|\|_{W_2^{1/2}(\mathbb
T)}\leq c_2\|f\|_{W_2^{1/2}(\mathbb T)}$ for all $f\in
W_2^{1/2}(\mathbb T)$, where $c_1, c_2>0$ do not depend on $f$.
Thus, we see that every function that satisfies the Lipschitz
condition of order greater then $1/2$ belongs to
$W_2^{1/2}(\mathbb T)$. We shall show that, in general, there
is no change of variable which will bring a complex-valued
function that satisfies the Lipschitz condition of order less
than $1/2$ into $W_2^{1/2}(\mathbb T)$. The author does not
know if the same holds for the functions satisfying the
Lipschitz condition of order $1/2$ (see Remarks at the end of
the paper).

\quad

\begin{center}
\textbf{2. Result}
\end{center}

Let $\omega$ be a modulus of continuity, i.e., a nondecreasing
continuous function on $[0, +\infty)$ such that $\omega(0)=0$
and $\omega(x+y)\leq\omega(x)+\omega(y)$. By
$\mathrm{Lip}_\omega(\mathbb T)$ we denote the class of all
complex-valued functions $f$ on $\mathbb T$ with $\omega(f,
\delta)=O(\omega(\delta)), ~\delta\rightarrow+0,$ where
$$
\omega(f, \delta)=\sup_{|t_1-t_2|\leq\delta} |f(t_1)-f(t_2)|,
\qquad \delta\geq 0,
$$
is the modulus of continuity of $f$. For $0<\alpha\leq 1$ we
just write $\mathrm{Lip}_\alpha$ instead of
$\mathrm{Lip}_{\delta^\alpha}$.

\quad

\textbf{Theorem.} \emph{Suppose that
$\limsup_{\delta\rightarrow+0}\omega(\delta)/\sqrt{\delta}=\infty$.
Then there exists a complex-valued function
$f\in\mathrm{Lip}_\omega(\mathbb T)$ such that $f\circ h\notin
W_2^{1/2}(\mathbb T)$ for every homeomorphism $h$ of the circle
$\mathbb T$ onto itself. In particular, if $\alpha<1/2$, then
there exists a function of class $\mathrm{Lip}_\alpha(\mathbb
T)$ with this property.}

\quad

Ideologically the method of the proof of this theorem is close
to the one used by the author to prove Theorem 4 in [7].

We shall need certain preliminary constructions and lemmas.
Simple Lemma 1 below is purely technical.

\quad

\textbf{Lemma 1.} \emph{Under the assumption of the theorem on
$\omega$ there exists a sequence $\delta_k>0, ~k=1, 2, \ldots,$
such that
$$
\sum_{k=1}^\infty\delta_k<2\pi/6,
\eqno(3)
$$
$$
\sum_{k=1}^\infty (\omega(\delta_k))^2=\infty.
\eqno(4)
$$
}

\quad

\emph{Proof.} For each $j=1, 2, \ldots$ we can find
$\varepsilon_j$ so that $0<\varepsilon_j<2^{-(j+1)}$ and
$$
\frac{(\omega(\varepsilon_j))^2}{\varepsilon_j}\geq 2^j.
$$
Chose positive integers $n_j$ satisfying
$$ \frac{1}{2^{j+1}\varepsilon_j}\leq
n_j<\frac{1}{2^j\varepsilon_j}, \qquad j=1, 2, \ldots.
$$
Let $N_0=1$ and let $N_j=N_{j-1}+n_j$ for $j=1, 2, \ldots$. We
define the sequence $\delta_k, ~k=1, 2, \ldots,$ by setting
$\delta_k=\varepsilon_j$ if $N_{j-1}\leq k<N_j, \,j=1, 2,
\ldots$. This yields
$$
\sum_{k=1}^\infty\delta_k=\sum_{j=1}^\infty\sum_{N_{j-1}\leq k< N_j}\delta_k=
\sum_{j=1}^\infty n_j\varepsilon_j\leq\sum_{j=1}^\infty\frac{1}{2^j}=1,
$$
and at the same time
$$
\sum_{N_{j-1}\leq k< N_j}(\omega(\delta_k))^2=n_j(\omega(\varepsilon_j))^2
\geq n_j\varepsilon_j2^j\geq\frac{1}{2}.
$$
The lemma is proved.

\quad

For a closed interval $I=[a, b]\subseteq (0, 2\pi)$ let
$\Delta_I$ denote the ``triangle'' function supported on $I$,
i.e., a continuous function on the interval $[0, 2\pi]$ such
that $\Delta_I(t)=0$ for all $t\in [0, a]\cup [b, 2\pi]$,
$\Delta_I(c)=1$, where $c=(a+b)/2$ is the center of $I$, and
$\Delta_I$ is linear on $[a, c]$ and on $[c, b]$.

Let $\delta_k, k=1, 2, \ldots,$ be the sequence from Lemma 1.
Consider intervals $I_k=[a_k, b_k]\subseteq (0,2\pi)$ of length
$b_k-a_k=6\delta_k$, where $a_k<b_k<a_{k+1}, ~k=1, 2, \ldots$
(see (3)). For each $k$ let $J_k$ denote the left half of
$I_k$, i.e., $J_k=[a_k, (a_k+b_k)/2], ~k=1, 2, \ldots$.

Everywhere below we use $u$ and $v$ to denote two real-valued
functions on $\mathbb T$ defined by
$$
u(t)=\sum_{k=1}^\infty \omega(\delta_k)\Delta_{I_k}(t),
\quad v(t)=\sum_{k=1}^\infty\omega(\delta_k)\Delta_{J_k}(t), \qquad t\in [0, 2\pi].
$$
We shall show that the function $f=u+iv$ satisfies the
assertion of the theorem.

\quad

\textbf{Lemma 2.} \emph{The functions $u$ and $v$ are of class
$\mathrm{Lip}_\omega(\mathbb T)$.}

\quad

\emph{Proof.} It is clear that for an arbitrary (closed)
interval $I\subseteq(0, 2\pi)$, the function $\Delta_I$
satisfies
$$
|\Delta_I(t_1)-\Delta_I(t_2)|\leq \frac{2}{|I|}|t_1-t_2|,
\qquad \textrm{for all}\quad t_1, t_2\in\mathbb [0, 2\pi],
\eqno(5)
$$
where $|I|$ is the length of $I$.

Note also that if $0<x\leq y$, then $\omega(y)/y\leq
2\omega(x)/x$. Indeed, let $n=[y/x]+1$, where $[\alpha]$
denotes the integer part of a number $\alpha$, then we have
$y\leq nx\leq 2y$, so
$$
\frac{\omega(y)}{y}\leq\frac{\omega(nx)}{y}\leq \frac{n\omega(x)}{y}\leq2\frac{\omega(x)}{x}.
$$

Let us show that $u\in\mathrm{Lip}_\omega(\mathbb T)$; for $v$
the proof is similar. It is easy to see that to prove the
inclusion $u\in\mathrm{Lip}_\omega(\mathbb T)$ it suffices to
verify that for all $t_1, t_2\in \bigcup_k I_k$ we have
$$
|u(t_1)-u(t_2)|\leq c\omega(|t_1-t_2|),
$$
where $c>0$ does not depend on $t_1$ and $t_2$.

First we consider the case when $t_1$ and $t_2$ belong to the
same interval $I_k$. If that is the case, then, since
$|t_1-t_2|\leq |I_k|=6\delta_k$, we have
$$
\frac{\omega(6\delta_k)}{6\delta_k}
\leq 2\frac{\omega(|t_1-t_2|)}{|t_1-t_2|},
$$
so (see (5)),
$$
|u(t_1)-u(t_2)|=\omega(\delta_k)|\Delta_{I_k}(t_1)-\Delta_{I_k}(t_2)|\leq
$$
$$
\leq\omega(\delta_k)\frac{2}{6\delta_k}|t_1-t_2|\leq
2\frac{\omega(6\delta_k)}{6\delta_k}|t_1-t_2|\leq 4\omega(|t_1-t_2|).
$$
Consider now the case when $t_1\in I_{k_1}, ~t_2\in I_{k_2},
~k_1\neq k_2$. We can assume that $t_1<t_2$, and hence
$0<t_1<b_{k_1}<a_{k_2}<t_2<2\pi$. Using the previous estimate,
we obtain
$$
|u(t_1)-u(t_2)|\leq |u(t_1)|+|u(t_2)|=
|u(t_1)-u(b_{k_1})|+|u(t_2)-u(a_{k_2})|\leq 8\omega(|t_1-t_2|).
$$
The lemma is proved.

\quad

For $n=1, 2, \ldots$ we define functions $u_n$ by
$$
u_n(t)=\max\{u(t), 1/n\}, \quad t\in\mathbb T.
$$

As above, $V(\mathbb T)$ stands for the class of functions of
bounded variation on $\mathbb T$.

\quad

\textbf{Lemma 3.} \emph{The functions $u_n, n=1, 2, \ldots,$
have the following properties:
$$
|u_n(t_1)-u_n(t_2)|\leq |u(t_1)-u(t_2)| \quad\textrm{for all}\quad t_1, t_2\in\mathbb T
\quad\textrm{and all}\quad n;
\eqno(6)
$$
$$
u_n\in V(\mathbb T)\quad \textrm{for all}\quad n;
\eqno(7)
$$
$$
\sup_n\bigg|\int_{\mathbb T}v(t) du_n(t)\bigg|=\infty.
\eqno(8)
$$}

\quad

\emph{Proof.} Properties (6) and (7) are obvious. Let us verify
(8). To this end consider the middle thirds of the intervals
$J_k$, namely, the intervals $J_k^*=[a_k+\delta_k,
a_k+2\delta_k], \,k=1, 2, \ldots$. Note that if
$$
\frac{\omega(\delta_k)}{3}\geq\frac{1}{n},
\eqno(9)
$$
then the function $u_n$ coincides with $u$ on $J_k^*$. So, if
(9) holds, then $u_n$ is monotonically increasing on
$J_k^\ast$, and for its values at the endpoints of $J_k^\ast$
we have
$$
u_n(a_k+\delta_k)=\omega(\delta_k)/3, \quad u_n(a_k+2\delta_k)=2\omega(\delta_k)/3.
$$
It is easily seen, that for each $k$
$$
\min_{J_k^*} v=2\omega(\delta_k)/3.
$$
Taking into account that $u$, and hence $u_n$, is
non-decreasing on each interval $J_k$, we see that for all $n$
and $k$ satisfying condition (9)
$$
\int_{J_k}v du_n\geq\int_{a_k+\delta_k}^{a_k+2\delta_k}v du_n\geq
\frac{2}{3}\omega(\delta_k)\int_{a_k+\delta_k}^{a_k+2\delta_k}du_n
=\frac{2}{3}\omega(\delta_k)\frac{1}{3}\omega(\delta_k)=
\frac{2}{9}(\omega(\delta_k))^2.
$$
In addition (since $u_n$ is non-decreasing on each $J_k$) we
have
$$
\int_{J_k}v du_n\geq 0
$$
for all $n$ and $k$. Thus, taking into account that $v$
vanishes outside $\bigcup_{k=1}^\infty J_k$, we obtain
$$
\int_{\mathbb T}v du_n=\sum_{k=1}^\infty\int_{J_k}v du_n\geq
\sum_{k \,:\, \omega(\delta_k)\geq 3/n}\int_{J_k}v du_n\geq
\sum_{k \,:\, \omega(\delta_k)\geq 3/n}\frac{2}{9}(\omega(\delta_k))^2.
$$
Applying (4) we see that (8) holds. The lemma is proved.

\quad

We shall also need the following auxiliary lemma.

\quad

\textbf{Lemma 4.} \emph{If $x, y\in W_2^{1/2}\cap C(\mathbb T)$
and $y\in V(\mathbb T)$, then
$$
\bigg|\frac{1}{2\pi}\int_\mathbb T x(t) dy(t)\bigg|\leq
\|x\|_{W_2^{1/2}(\mathbb T)}\|y\|_{W_2^{1/2}(\mathbb T)}.
$$}

\quad

\emph{Proof.} Integration by parts yields
$$
\frac{1}{2\pi}\int_0^{2\pi} e^{ikt} dy(t)=-\frac{1}{2\pi}\int_0^{2\pi} y(t)de^{ikt}=
-ik\widehat{y}(-k).
$$
So, if $x$ is a trigonometric polynomial, then, using Cauchy
inequality, we obtain
$$
\bigg|\frac{1}{2\pi}\int_\mathbb T x(t) dy(t)\bigg|=
\bigg|\sum_k\widehat{x}(k)\frac{1}{2\pi}\int_\mathbb T e^{ikt} dy(t)\bigg|=
$$
$$
=\bigg|\sum_{k}\widehat{x}(k)(-ik)\widehat{y}(-k)\bigg|
\leq\|x\|_{W_2^{1/2}(\mathbb T)}\|y\|_{W_2^{1/2}(\mathbb T)}.
$$
To see that the assertion of the lemma holds in the general
case, consider the Fej\'er sums $\sigma_N(x)$ of the function
$x$:
$$
\sigma_N(x)(t)=\sum_{|k|\leq N}\bigg(1-\frac{|k|}{N}\bigg)\widehat{x}(k)e^{ikt}.
$$
Since $|\widehat{\sigma_N(x)}(k)|\leq|\widehat{x}(k)|$ for all
$k\in\mathbb Z$, we have $\|\sigma_N(x)\|_{W_2^{1/2}(\mathbb
T)}\leq\|x\|_{W_2^{1/2}(\mathbb T)}$. Hence,
$$
\bigg|\frac{1}{2\pi}\int_\mathbb T \sigma_N(x)(t) dy(t)\bigg|\leq
\|\sigma_N(x)\|_{W_2^{1/2}(\mathbb T)}\|y\|_{W_2^{1/2}(\mathbb T)}\leq
\|x\|_{W_2^{1/2}(\mathbb T)}\|y\|_{W_2^{1/2}(\mathbb T)}.
$$
At the same time, since $y$ is of bounded variation and
$\sigma_N(x)$ converges uniformly to $x$ it is clear that
$$
\frac{1}{2\pi}\int_\mathbb T\sigma_N(x)(t)dy(t)\rightarrow\frac{1}{2\pi}\int_\mathbb T x(t)dy(t)
$$
as $N\rightarrow\infty$. The lemma is proved.

\quad

Now we proceed directly to the proof of the theorem. Let
$f=u+iv$. Lemma 2 yields $f\in\mathrm{Lip}_\omega (\mathbb T)$,
so it remains to show that $f\circ h\notin W_2^{1/2}(\mathbb
T)$ for every homeomorphism $h$ of $\mathbb T$. It is obvious
that if a function is in $W_2^{1/2}(\mathbb T)$, then both its
real and imaginary parts are in $W_2^{1/2}(\mathbb T)$ as well.
Assume that, contrary to the assertion of the theorem, $f\circ
h\in W_2^{1/2}(\mathbb T)$ for a certain homeomorphism $h$.
Then we have $u\circ h\in W_2^{1/2}(\mathbb T)$ and $v\circ
h\in W_2^{1/2}(\mathbb T)$.

Note that (6) implies $|u_n\circ h(t_1)-u_n\circ h(t_2)|\leq
|u\circ h(t_1)-u\circ h(t_2)|$ for all $t_1, t_2\in\mathbb T$.
Using the equivalence of the semi-norms
$\|\cdot\|_{W_2^{1/2}(\mathbb T)}$ and
$|\|\cdot|\|_{W_2^{1/2}(\mathbb T)}$ (see (2)), we infer that
$u_n\circ h\in W_2^{1/2}(\mathbb T)$ for all $n=1, 2, \ldots$,
and
$$
\|u_n\circ h\|_{W_2^{1/2}(\mathbb T)}\leq
c \|u\circ h\|_{W_2^{1/2}(\mathbb T)}, \qquad n=1, 2, \ldots,
\eqno(10)
$$
where $c>0$ does not depend on $n$.

The property of a function to be of bounded variation is
invariant under homeomorphic changes of variable, hence from
(7) it follows that $u_n\circ h\in V(\mathbb T)$ for all $n$.
Certainly we also have $u_n\circ h\in C(\mathbb T)$. Applying
Lemma 4, and taking (10) into account, we obtain
$$
\bigg|\frac{1}{2\pi}\int_\mathbb T v(t) du_n(t)\bigg|=
\bigg|\frac{1}{2\pi}\int_\mathbb T v\circ h(t) du_n\circ h(t)\bigg|\leq
$$
$$
\leq\|v\circ h\|_{W_2^{1/2}(\mathbb T)} \|u_n\circ
h\|_{W_2^{1/2}(\mathbb T)}\leq c \|v\circ h\|_{W_2^{1/2}(\mathbb T)} \|u\circ
h\|_{W_2^{1/2}(\mathbb T)},
$$
which contradicts (8). The theorem is proved.

\quad

\emph{Remarks.} 1. For $s>0$ consider the Sobolev space
$W_2^s(\mathbb T)$ i.e., the space of all (integrable)
functions $f$ with
$$
\sum_{k\in\mathbb Z}|\widehat{f}(k)|^2|k|^{2s}<\infty.
$$
As the author of this paper showed in [7, Corollary 3], for
each compact family $K$ in $C(\mathbb T)$ (or equivalently for
each class $\mathrm{Lip}_\omega(\mathbb T)$) there exists a
homeomorphism $h$ of $\mathbb T$ such that $f\circ h\in
\bigcap_{s<1/2}W_2^s(\mathbb T)$ for all $f\in K$ (for all
$f\in\mathrm{Lip}_\omega(\mathbb T)$).

2. There exists a real-valued $f\in C(\mathbb T)$ such that
$f\circ h\notin\bigcup_{s>1/2}W_2^s(\mathbb T)$ for every
homeomorphism $h$ of $\mathbb T$. This is a simple consequence
of the inclusion $\bigcup_{s>1/2}W_2^s\cap C(\mathbb
T)\subseteq A(\mathbb T)$, where $A(\mathbb T)$ is the Wiener
algebra of absolutely convergent Fourier series, and a
well-known result of Olevski\v{\i}, that provides a negative
answer to Lusin's rearrangement problem: there exists a
real-valued $f\in C(\mathbb T)$ such that $f\circ h\notin
A(\mathbb T)$ for every homeomorphism $h$ ([8], see also [9]).

3. The function $f(t)=\sum_{k\geq 0}2^{-k/2}e^{i2^kt}$ is in
$\mathrm{Lip}_{1/2}(\mathbb T)$, (see, e.g., [1, Ch. XI, Sec.
6]), but it is obvious, that $f\notin W_2^{1/2}(\mathbb T)$;
thus $\mathrm{Lip}_{1/2}(\mathbb T)\nsubseteq W_2^{1/2}(\mathbb
T)$. The author does not know if the assertion of the theorem
proved in this paper holds for $\omega(\delta)=\delta^{1/2}$.
At the same time there is no change of variable which will
bring the whole class $\mathrm{Lip}_{1/2}(\mathbb T)$ into
$W_2^{1/2}(\mathbb T)$; a proof will be presented in another
paper.

\quad

\quad

\quad

\begin{center}
\textbf{References}
\end{center}

\flushleft
\begin{enumerate}

\item Bari N. K., \emph{Trigonometric series}, Fizmatgiz,
    Moscow, 1961. English transl.: Vols. I, II, Pergamon
    Press, Oxford, and Macmillan, New York, 1964.

\item Bohr H., \"Uber einen Satz von J. P\'al, \emph{Acta
    Sci. Math. (Szeged)} \textbf{7} (1935), 129--135.

\item Goffman G., Nishiura T., Waterman D.,
    \emph{Homeomorphisms in Analysis}, Mathematical Surveys
    and Monographs v. 54, Amer. Math. Soc., 1997.

\item Jurkat W., Waterman D., Conjugate functions and the
    Bohr--P\'al theorem, \emph{Complex Variables} \textbf{12}
    (1989), 67--70.

\item Kahane J.-P., Quatre le\c cons sur les
    hom\'eomorphismes du circle et les s\'eries de Fourier,
    in: \emph{Topics in Modern Harmonic Analysis}, Vol. II,
    Ist. Naz. Alta Mat. Francesco Severi, Roma, 1983,
    955--990.

\item Kahane J.-P., Katznelson Y., S\'eries de Fourier des
    fonctions born\'ee, \emph{Studies in Pure Math.}, in
    Memory of Paul Tur\'an, Budapest, 1983, pp. 395--410.
    (Preprint, Orsay, 1978.)

\item Lebedev  V. V., Change of variable and the rapidity of
    decrease of Fourier coefficients, \emph{Mat. Sbornik}
    \textbf{181} (1990), 1099--1113 (Russian). English
    transl.: Math. USSR-Sb. \textbf{70} (1991), 541--555.
    English transl. corrected by the author is available at:
    http://arxiv.org/abs/1508.06673

\item Olevski\v{\i} A. M., Change of variable and absolute
    convergence of Fourier series, \emph{Dokl. Akad. Nauk
    SSSR} \textbf{256} (1981), 284--288 (Russian). English
    transl.: \emph{Soviet Math. Dokl.} \textbf{23} (1981),
    76--79.

\item Olevski\v{\i} A. M., Modifications of functions and
    Fourier series, \emph{Uspekhi Mat. Nauk} \textbf{40}
    (1985), 157--193 (Russian). English transl.:
    \emph{Russian Math. Surveys} \textbf{40} (1985),
    181--224.

\item Olevski\v{\i} A. M., Modifications of functions and
    Fourier series, \emph{Theory of Functions and
    Approximations}, Proc. Second Saratov Winter School
    (Saratov, Jan 24--Feb 5, 1984), Part 1, Saratov. Gos.
    Univ., Saratov, 1986, pp. 31--43 (Russian).

\item P\'al J., Sur les transformations des fonctions qui
    font converger leurs s\'eries de Fourier, \emph{C. R.
    Acad. Sci Paris} \textbf{158} (1914), 101--103.

\item Saakjan A. A., Integral moduli of smoothness and the
    Fourier coefficients of the composition of functions,
    \emph{Mat. Sbornik} \textbf{110} (1979), 597--608
    (Russian). English transl.: \emph{Mathematics of the
    USSR-Sbornik} \textbf{38} (1981), 549--561.
	
\end{enumerate}

\quad

Moscow Institute of Electronics and Mathematics,

National Research University Higher School of Economics

34 Tallinskaya St.

Moscow, 123458, Russia

E-mail address: \emph{lebedevhome@gmail.com}

\end{document}